\documentclass[12pt,a4paper,fleqn,leqno]{amsart}
\usepackage{amssymb}
\usepackage{mathrsfs}

\theoremstyle{plain}
\newtheorem{theorem}{Theorem}[section]

\newtheorem{proposition}[theorem]{Proposition}
\newtheorem{claim}[theorem]{Claim}

\theoremstyle{definition}
\newtheorem{problem}{{Question}}

\renewcommand{\emptyset}{\varnothing}

\DeclareMathOperator{\uhr}{\upharpoonright}

\DeclareMathOperator{\setR}{\mathbb{R}}

\DeclareMathOperator{\diam}{diam}
\DeclareMathOperator{\conv}{conv}

\numberwithin{equation}{section} 

\begin{document}

\title{Selections, Extensions and Collectionwise Normality}

\author{Valentin Gutev}

\address{School of Mathematical Sciences, University of KwaZulu-Natal,
  Westville Campus, Private Bag X54001, Durban 4000, South Africa}

\email{gutev@ukzn.ac.za}

\thanks{The work of the first author is based upon research supported
  by the NRF of South Africa}

\author{Narcisse Roland Loufouma Makala}

\address{School of Mathematical Sciences, University of KwaZulu-Natal,
  Westville Campus, Private Bag X54001, Durban 4000, South Africa}

\email{roland@aims.ac.za}



\begin{abstract}
  We demonstrate that the classical Michael selection theorem for
  l.s.c.\ mappings with a collectionwise normal domain can be reduced
  only to compact-valued mappings modulo Dowker's extension theorem
  for such spaces. The idea used to achieve this reduction is also
  applied to get a simple direct proof of that selection theorem
  of Michael's. Some other possible applications are demonstrated as
  well.
\end{abstract}

\date{August 18, 2009}
\maketitle

\section{Introduction}
\label{section-introduction}

For a topological space $E$, let $2^E$ be the family of all nonempty
subsets of $E$; $\mathscr{F}(E)$ --- the subfamily of $2^E$ consisting
of all closed members of $2^E$; and let $\mathscr{C}(E)$ be that one of
all compact members of $\mathscr{F}(E)$. Also, let
$\mathscr{C}'(E)=\mathscr{C}(E)\cup \{E\}$.  A set-valued mapping $\varphi:X\to
2^E$ is \emph{lower semi-continuous}, or l.s.c., if the set
\[
\varphi^{-1}(U)=\{x\in X: \varphi(x)\cap U\neq \emptyset\}
\]
is open in $X$ for every open $U\subset E$. A map $f:X\to E$ is a
\emph{selection} for $\varphi:X\to 2^E$ if $f(x)\in \varphi(x)$ for every $x\in
X$. \medskip

Recall that a space $X$ is \emph{collectionwise normal} if it is a
$T_1$-space and for every discrete collection $\mathscr{D}$ of closed
subsets of $X$ there exists an open discrete family $\big\{U_D: D\in
\mathscr{D}\big\}$ such that $D\subset U_D$ for every $D\in \mathscr{D}$. Every
collectionwise normal space is normal, but the converse is not
necessarily true \cite{bing:51}, see, also, \cite[5.1.23 Bing's
Example]{engelking:89}. It is well known that a $T_1$-space $X$ is
collectionwise normal if and only if for every closed subset $A\subset X$,
every continuous map from $A$ to a Banach space $E$ can be
continuously extended to the whole of $X$, Dowker
\cite{dowker:52}. Generalizing this result, Michael \cite{michael:56a}
stated the following theorem.

\begin{theorem}[\cite{michael:56a}]
  \label{theorem-michael}
  For a $T_1$-space $X$, the following are equivalent\textup{:}
  \begin{enumerate}
  \item[(a)] $X$ is collectionwise normal.
  \item[(b)] If $E$ is a Banach space and $\varphi:X\to \mathscr{C}'(E)$ is an
    l.s.c.\ convex-valued mapping, then $\varphi$ has a continuous
    selection.
  \end{enumerate}
\end{theorem}

However, the arguments in \cite{michael:56a} for (a)$\ \Rightarrow\ $(b) of
Theorem \ref{theorem-michael} were incomplete and, in fact, working
only for the case of $\mathscr{C}(E)$-valued mappings. The first
complete proof of this implication was given by Choban and Valov
\cite{choban-valov:75} using a different technique.  We are now ready
to state also the main purpose of this paper. Namely, in this paper we
prove the following theorem which demonstrates that the original
Michael arguments in \cite{michael:56a} have been actually adequate to
the proof of Theorem~\ref{theorem-michael}.

\begin{theorem}
  \label{theorem-main-1}
  For a Banach space $E$, the following are equivalent\textup{:}
  \begin{enumerate}
  \item[(a)] If $X$ is a collectionwise normal space and ${\varphi:X\to
      \mathscr{C}(E)}$ is an l.s.c.\ convex-valued mapping, then $\varphi$
    has a continuous selection.
  \item[(b)] If $X$ is a collectionwise normal space and $\varphi:X\to
    \mathscr{C}'(E)$ is an l.s.c.\ convex-valued mapping, then $\varphi$ has
    a continuous selection.
  \end{enumerate}
\end{theorem}

Let us emphasize that the proof of Theorem \ref{theorem-main-1} is
based only on Dowker's extension theorem \cite{dowker:52}. This proof
is presented in the next section and its main ingredient is the fact
that if $\varphi:X\to \mathscr{C}'(E)$ and $g:X\to E$, then $\varphi(x)$ is compact
for every $x\in X$ for which $g(x)\notin \varphi(x)$. This is further applied in
Section \ref{sec:more-select-probl} to get with ease a direct proof of
a natural generalization of
Theorem~\ref{theorem-michael}. Section~\ref{sec:contr-select-count}
deals with controlled selections for set-valued mappings defined on
countably paracompact or collectionwise normal spaces which are
naturally interrelated to the idea of Theorem \ref{theorem-main-1}.

\section{Proof of Theorem \ref{theorem-main-1}}
\label{section-proof-main-theorem}

It only suffices to prove that (a)$\ \Rightarrow\ $(b). To this end, suppose
that (a) of Theorem \ref{theorem-main-1} holds, $E$ is a Banach space,
and $X$ is a collectionwise normal space. Here, and in the sequel, we
will use $d$ to denote the metric on $E$ generated by the norm of
$E$. Recall that a map $f:X\to E$ is an \emph{$\varepsilon$-selection} for $\psi :X\to
2^E$ if $d(f(x),\psi (x))<\varepsilon$ for every $x\in X$. \medskip

The key element in the proof of this implication is the following
construction of approximate selections.

\begin{claim}
  \label{claim-main-construction}
  Let $\psi:X\to 2^E$ be an l.s.c.\ convex-valued mapping and $g:X\to E$ be a
  continuous map such that $\psi(x)$ is compact whenever $x\in X$ and
  $g(x)\notin \psi(x)$. Then, for every $\varepsilon>0$, $\psi$ has a continuous
  $\varepsilon$-selection.
\end{claim}

\begin{proof}
  Let $\varepsilon>0$ and $A=\big\{x\in X: d(g(x),\psi(x))\geq\varepsilon \big\}$. Then, $A\subset X$ is
  closed because $\psi$ is l.s.c.\ and $g$ is continuous. Since $\psi \uhr
  A:A\to \mathscr{C}(E)$ and $A$ is itself a collectionwise normal
  space, by (a) of Theorem \ref{theorem-main-1}, $\psi \uhr A$ has a
  continuous selection $h_0:A\to E$. Since $X$ is collectionwise normal,
  by Dowker's extension theorem \cite{dowker:52}, there exists a
  continuous map $h:X\to E$ such that $h\uhr A=h_0$. Consider the set $
  U=\big\{x\in X: d(h(x),\psi(x))<\varepsilon \big\}$ which contains $A$ and is open
  because $\psi$ is l.s.c.\ and $h$ is continuous. Finally, take a
  continuous function $\alpha:X\to [0,1]$ such that $A\subset \alpha^{-1}(0)$ and $X\setminus U\subset
  \alpha^{-1}(1)$, and then define a continuous map $f:X\to E$ by
  \[
  f(x)=\alpha(x)\cdot g(x)+\big(1-\alpha(x)\big)\cdot h(x),\quad x\in X.
  \]
  This $f$ is as required. Indeed, take a point $x\in X$. If $x\in A$,
  then $\alpha(x)=0$ and, therefore, $f(x)=h(x)=h_0(x)\in \psi(x)$. If $x\in X\setminus
  U$, then $\alpha(x)=1$ and we now have that $f(x)=g(x)$, so
  $d(f(x),\psi(x))=d(g(x),\psi(x))<\varepsilon$ because $x\notin A$. Suppose finally that
  $x\in U\setminus A$. In this case, $d(h(x),\psi(x))<\varepsilon$ and
  $d(g(x),\psi(x))<\varepsilon$. Since $\psi(x)$ is convex and $f(x)=\alpha(x)\cdot
  g(x)+\big(1-\alpha(x)\big)\cdot h(x)$, this implies that
  $d(f(x),\psi(x))<\varepsilon$. The proof is completed.
\end{proof}

Having already established Claim \ref{claim-main-construction}, we
proceed to the proof of (a)$\ \Rightarrow\ $(b) which is based on standard
arguments for constructing continuous selections, see
\cite{michael:56a}. In this proof, and in what follows, for a nonempty
subset $S\subset E$ and $\varepsilon>0$, we will use $B_\varepsilon(S)=\{y\in E: d(y,S)<\varepsilon\}$ to
denote the \emph{open $\varepsilon$-neighbourhood} of $S$. In particular, for a
point $y\in E$, we set $B_\varepsilon(y)=B_\varepsilon(\{y\})$.\smallskip

Let $\varphi:X\to \mathscr{C}'(E)$ be an l.s.c.\ convex-valued mapping. If
$g:X\to E$ is any continuous map, say a constant one, then $\varphi(x)$ is
compact for every $x\in X$ for which $g(x)\notin \varphi(x)$. Hence, by Claim
\ref{claim-main-construction}, $\varphi$ has a continuous $2^{-1}$-selection
$f_0:X\to E$. Define $\varphi_1:X\to \mathscr{F}(E)$ by
\[
\varphi_1(x)=\overline{\varphi(x)\cap B_{2^{-1}}(f_0(x))},\quad x\in X.
\]
According to \cite[Propositions 2.3 and 2.5]{michael:56a}, $\varphi_1$ is
l.s.c., and clearly it is convex-valued. Finally, observe that if
$f_0(x)\notin \varphi_1(x)$ for some $x\in X$, then $f_0(x)\notin \varphi(x)$ and, therefore,
$\varphi(x)$ is compact because it is $\mathscr{C}'(E)$-valued. Since
$\varphi_1(x)$ is a closed subset of $\varphi(x)$, it is also compact. Hence, by
Claim \ref{claim-main-construction}, $\varphi_1$ has a continuous
$2^{-2}$-selection $f_1$. In particular, $f_1$ is a continuous
$2^{-2}$-selection for $\varphi$ such that
\[
d(f_1(x),f_0(x))\leq 2^{-1}<2^0,\quad \text{for every $x\in X$.}
\]
Thus, by induction, we get a sequence $\{f_n:n<\omega\}$ of continuous maps
such that, for every $n<\omega$ and $x\in X$,
\begin{align}
  \label{eq:fundamental-1}
  d(f_n(x),\varphi(x))&<2^{-(n+1)},\\
  \label{eq:fundamental-2}
  d(f_{n+1}(x),f_{n}(x))&<2^{-n}.
\end{align}
By \eqref{eq:fundamental-2}, $\{f_n:n<\omega\}$ is a Cauchy sequence in the
complete metric space $(E,d)$, so it must converge to some continuous
$f:X\to E$. By \eqref{eq:fundamental-1}, $f(x)\in \varphi(x)$ for every $x\in
X$. Hence, (b) holds and the proof of Theorem \ref{theorem-main-1} is
completed.

\section{More on Selections and Collectionwise Normality}
\label{sec:more-select-probl}

A space $X$ is \emph{$\tau$-collectionwise normal}, where $\tau$ is an
infinite cardinal number, if it is a $T_1$-space and for every
discrete collection $\mathscr{D}$ of closed subsets of $X$, with
$|\mathscr{D}|\leq \tau$, there exists a discrete collection $\{U_D : D\in
\mathscr{D}\}$ of open subsets of $X$ such that $D\subset U_D$ for every $D\in
\mathscr{D}$.  Clearly, a space $X$ is collectionwise normal if and
only if it is $\tau$-collectionwise normal for every $\tau$. Also, it is
well known that $X$ is normal if and only if it is $\omega$-collectionwise
normal. However, for every $\tau$ there exists a $\tau$-collectionwise
normal space which is not $\tau^+$-collectionwise normal
\cite{przymusinski:78a}, where $\tau^+$ is the immediate successor of
$\tau$.\medskip

The proof of Theorem \ref{theorem-main-1} suggests an easy direct
proof of the following natural generalization of Theorem
\ref{theorem-michael} in \cite{choban-valov:75}.

\begin{theorem}[\cite{choban-valov:75}]
  \label{theorem-choban-valov-direct}
  Let $X$ be a $\tau$-collectionwise normal space, $E$ be a Banach space
  with a topological weight $w(E)\leq \tau$, and let $\varphi:X\to \mathscr{C}'(E)$
  be an l.s.c.\ convex-valued mapping. Then, $\varphi$ has a continuous
  selection.
\end{theorem}

\begin{proof}
  It only suffices to prove the statement of Claim
  \ref{claim-main-construction} for this particular case. So, suppose
  that $\psi:X\to 2^E$ is l.s.c.\ and convex-valued, and $g:X\to E$ is a
  continuous map such that $\psi(x)$ is compact whenever $g(x)\notin
  \psi(x)$. Also, let $\varepsilon>0$ and let $\mathscr{V}$ be an open and locally
  finite cover of $E$ such that $\diam_d(V)<\varepsilon$ for every $V\in
  \mathscr{V}$. Since $g$ is continuous, $\mathscr{U}_1=\{g^{-1}(V)\cap
  \psi^{-1}(V): V\in \mathscr{V}\}$ is a locally finite family of open
  subsets of $X$ which refines $\{\psi^{-1}(V): V\in \mathscr{V}\}$. Then,
  $A=X\setminus \bigcup \mathscr{U}_1$ is a closed subset of $X$, while $\psi\uhr A$ is
  compact-valued. Indeed, if $g(x)\in \psi(x)$, then $x\in \psi^{-1}(V)$
  whenever $V\in \mathscr{V}$ and $g(x)\in V$. That is, $x\in A$ implies
  $g(x)\notin \psi(x)$, so, in this case, $\psi(x)$ must be compact. Thus,
  $\{\psi^{-1}(V): V\in \mathscr{V}\}$ is an open (in $X$) and point-finite
  (in $A$) cover of $A$ such that $|\mathscr{V}|\leq \tau$ because
  $\mathscr{V}$ is locally-finite and $w(E)\leq \tau$. Since $X$ is
  $\tau$-collectionwise normal, by \cite[Lemma 1.6]{nedev:80},
  $\{\psi^{-1}(V): V\in \mathscr{V}\}$ has an open and locally finite (in
  $X$) refinement $\mathscr{U}_2$ which covers $A$. Then,
  $\mathscr{U}=\mathscr{U}_1\cup \mathscr{U}_2$ is an open and locally
  finite cover of $X$ which refines $\{\psi^{-1}(V): V\in \mathscr{V}\}$. For
  every $U\in \mathscr{U}$ take a fixed $V_U\in \mathscr{V}$ such that $U\subset
  \psi^{-1}(V_U)$ and a point $e(U)\in V_U$ provided $V_U\neq \emptyset$. Next, take a
  partition of unity $\{\xi_U: U\in \mathscr{U}\}$ on $X$ which is index
  subordinated to the cover $\mathscr{U}$, see
  \cite{michael:53}. Finally, define a continuous map $f:X\to E$ by
  $f(x)=\sum\{\xi_U(x)\cdot e(U): U\in \mathscr{U}\}$, $x\in X$. This $f$ is an
  $\varepsilon$-selection for $\psi$.
\end{proof}

As far as the role of the family $\mathscr{C}'(E)$ is concerned, the
arguments in the proof of Theorems \ref{theorem-main-1} and
\ref{theorem-choban-valov-direct} were based only on the property that
if $\varphi:X\to \mathscr{C}'(E)$ and $g:X\to E$ is an $\varepsilon$-selection for $\varphi$ for
some $\varepsilon>0$, then the set-valued mapping $\psi(x)=\overline{\varphi(x)\cap
  B_\varepsilon(g(x))}$, $x\in X$, is such that $\psi(x)$ is compact whenever $g(x)\notin
\psi(x)$. That is, this resulting $\psi$ is always as in Claim
\ref{claim-main-construction}, and the inductive construction can be
carried on.\medskip

Motivated by this, we shall say that a mapping $\psi:X\to \mathscr{F}(E)$
has a \emph{selection $\mathscr{C}(E)$-deficiency} if there exists a
continuous $g:X\to E$ such that $\psi(x)\in \mathscr{C}(E)$ for every $x\in X$
for which $g(x)\notin \psi(x)$. Clearly, every $\varphi:X\to \mathscr{C}'(E)$ has this
property, for instance take $g:X\to E$ to be any constant map. However,
there are natural examples of mappings $\varphi:X\to \mathscr{F}(E)$ which
have a selection $\mathscr{C}(E)$-deficiency and are not
$\mathscr{C}'(E)$-valued, see next section.  Related to this, we don't
know if such mappings may have continuous selections in the case of
collectionwise normal spaces.

\begin{problem}
  \label{question-selection-deficiency}
  Let $X$ be a collectionwise normal space, $E$ be a Banach space, and
  let $\varphi:X\to \mathscr{F}(E)$ be an l.s.c.\ convex-valued mapping which
  has a selection $\mathscr{C}(E)$-deficiency. Then, is it true that
  $\varphi$ has a continuous selection?
\end{problem}

Another aspect of improving Theorem \ref{theorem-michael} is related
to the range of the set-valued mapping. In this regard, Theorem
\ref{theorem-choban-valov-direct} remains valid without any change in
the arguments if the Banach space $E$ is replaced by a closed convex
subset $Y$ of $E$. On the other hand, if $Y$ is a completely
metrizable absolute retract for the metrizable spaces, then for every
collectionwise normal space $X$ and closed $A\subset X$, every continuous
map $g:A\to Y$ can be continuously extended to the whole of $X$, see,
e.g., \cite{przymusinski:78a}. In particular, this is true for every
convex $G_\delta$-subset $Y$ of a Banach space $E$. Namely, $Y$ is an
absolute retract for metrizable spaces being convex (by Dugundji's
extension theorem \cite{dugundji:51}), and is also completely
metrizable being a $G_\delta$-subset of a complete metric space. Motivated
by this and the relations between extensions and selections
demonstrated in the proof of Theorem \ref{theorem-main-1}, we have
also the following question.

\begin{problem}
  \label{question-michael-related}
  Let $E$ be a Banach space, $Y\subset E$ be a convex $G_\delta$-subset of $E$,
  $X$ be a collectionwise normal space, and let $\varphi:X\to \mathscr{C}'(Y)$
  be an l.s.c.\ convex-valued mapping. Then, is it true that $\varphi$ has a
  continuous selection?
\end{problem}

Question \ref{question-michael-related} is similar to Michael's
$G_\delta$-problem \cite[Problem 396]{michael:90} if for a paracompact
space $X$ and a convex $G_\delta$-subset $Y$ of a Banach space, every
l.s.c.\ convex-valued $\varphi:X\to \mathscr{F}(Y)$ has a continuous
selection. In general, the answer to this latter problem is in the
negative due to a counterexample constructed by Filippov
\cite{filippov:04,filippov:05}. However, Michael's $G_\delta$-problem was
resolved in the affirmative in a number of partial cases. The solution
in some of these cases remains valid for Question
\ref{question-michael-related} as well. For instance, if $Y$ is a
countable intersection of open convex sets, then the closure
convex-hull $\overline{\conv(K)}$ of every compact subset of $Y$ will
be still a subset of $Y$, see \cite{michael:90}. In this case, by a
result of \cite{choban-valov:75}, $\varphi$ will have an l.s.c.\
convex-valued selection $\psi:X\to \mathscr{C}(Y)$ (i.e., $\psi(x)\subset \varphi(x)$ for
all $x\in X$). Hence, $\varphi$ will have a continuous selection because, by
Theorem \ref{theorem-michael}, so does $\psi$. If the covering dimension
of $X$ is bounded (i.e., $\dim(X)<\infty$), then the answer to Question
\ref{question-michael-related} is also ``yes'', this follows directly
from a selection theorem in \cite{gutev:86}. The answer to Question
\ref{question-michael-related} is also ``yes'' if $X$ is strongly
countable-dimensional (i.e., a countable union of closed
finite-dimensional subsets). In this case, there exists a metrizable
(strongly) countable-dimensional space $Z$, a continuous map $g:X\to Z$
and an l.s.c.\ mapping $\psi:Z\to \mathscr{C}(Y)$ such that $\psi(g(x))\subset
\varphi(x)$ for every $x\in X$, see, for instance, the proof of \cite[Theorem
5.3]{nedev:80}. Next, define a mapping $\Phi:Z\to \mathscr{F}(Y)$ by
$\Phi(z)=\overline{\conv(\psi(z))}^Y$, $z\in Z$, where the closure is in
$Y$. According to \cite[Propositions 2.3 and 2.6]{michael:56a}, $\Phi$
remains l.s.c., and, by \cite[Corollary 1.2]{gutev:94}, it admits a
continuous selection $h:Z\to Y$. Then, $f=h\circ g$ is a continuous
selection for $\varphi$ because $\Phi(g(x))\subset \varphi(x)$ for all $x\in X$.

\section{Controlled Selections and Countable Paracompactness}
\label{sec:contr-select-count}

A function $\xi:X\to \setR$ is \emph{lower} (\emph{upper})
\emph{semi-continuous} if the set
\[
\{x\in X:\xi(x)>r\}\quad \text{(respectively, $\{x\in X:\xi(x)<r\}$)}
\]
is open in $X$ for every $r\in \setR$. If $(E,d)$ is a metric space, $\varphi:X\to
2^E$ and $\eta:X\to (0,+\infty)$, then we shall say that $g:X\to E$ is an
\emph{$\eta$-selection} for $\varphi$ if $d(g(x),\varphi(x))<\eta(x)$ for every $x\in
X$. \medskip

In this section, we first prove the following characterization of
countably paracompact normal spaces.

\begin{theorem}
  \label{theorem-countably-paracompact-controlled}
  For a $T_1$-space $X$, the following are equivalent\textup{:}
  \begin{enumerate}
  \item[(a)] $X$ is countably paracompact and normal.
  \item[(b)] If $E$ is a separable Banach space, $\varphi:X\to \mathscr{F}(E)$
    is an l.s.c.\ convex-valued mapping, $\eta:X\to (0,+\infty)$ is lower
    semi-continuous and $g:X\to E$ is a continuous $\eta$-selection for
    $\varphi$, then $\varphi$ has a continuous selection $f:X\to E$ such that
    $d(f(x),g(x))<\eta(x)$ for all $x\in X$.
  \item[(c)] If $\varphi:X\to \mathscr{C}(\setR)$ is an l.s.c.\ convex-valued
    mapping, $\varepsilon>0$ and $g:X\to \setR$ is a continuous $\varepsilon$-selection for $\varphi$,
    then $\varphi$ has a continuous selection $f:X\to \setR$ such that
    $d(f(x),g(x))<\varepsilon$ for all $x\in X$.
  \end{enumerate}
\end{theorem}

\begin{proof}
  (a)$\ \Rightarrow\ $(b). Let $X$ be a countably paracompact normal space, and
  let $E$, $\varphi$, $\eta$ and $g$ be as in (b). Since $\varphi$ is l.s.c.\ and $g$
  is continuous, $\xi(x)=d(g(x),\varphi(x))$, $x\in X$, is an upper
  semi-continuous function such that $\xi(x)<\eta(x)$ for all $x\in X$
  because $g$ is an $\eta$-selection for $\varphi$. Since $X$ is countably
  paracompact and normal, by a result of
  \cite{dieudonne:44,dowker:51,katetov:51} (see, also,
  \cite[5.5.20]{engelking:89}) there exists a continuous function
  $\alpha:X\to \setR$ such that $\xi(x)<\alpha(x)<\eta(x)$ for every $x\in X$. Then, define
  an l.s.c.\ mapping $\psi:X\to \mathscr{F}(E)$ by $\psi(x)=\overline{\varphi(x)\cap
    B_{\alpha(x)}(g(x))}$, $x\in X$. Since $\psi$ is convex-valued, by
  \cite[Theorem 3.1$''$]{michael:56a}, $\psi$ has a continuous selection
  $f:X\to E$. In particular, $d(f(x),g(x))\leq \alpha(x)<\eta(x)$ for all $x\in
  X$. \smallskip

  Since (b)$\ \Rightarrow\ $(c) is obvious, we complete the proof showing that
  (c)$\ \Rightarrow\ $(a). So, suppose that (c) holds. If $A$ and $B$ are
  disjoint closed subsets of $X$, then $\varphi(x)=\{0\}$ if $x\in A$,
  $\varphi(x)=\{1\}$ if $x\in B$, and $\varphi(x)=[0,1]$ otherwise, is an l.s.c.\
  convex-valued mapping $\varphi:X\to \mathscr{C}(\setR)$. If $g(x)=\frac 12$, $x\in
  X$, then $g$ is a continuous $1$-selection for $\varphi$, and, by (c), $\varphi$
  has a continuous selection $f:X\to \setR$. According to the definition of
  $\varphi$, we get that $A\subset f^{-1}(0)$ and $B\subset f^{-1}(1)$, hence $X$ is
  normal. To show that $X$ is countably paracompact, let $\{F_n:n<\omega\}$
  be a decreasing sequence of closed subsets of $X$ such that $F_0=X$
  and $\bigcap\{F_n:n<\omega\}=\emptyset$. Next, for every $x\in X$, let $n(x)=\max\{n<\omega: x\in
  F_n\}$. Then, define a convex-valued mapping $\varphi:X\to \mathscr{C}(\setR)$ by
  $\varphi(x)=\left[0,2^{-n(x)}\right]$, $x\in X$. Observe that $\varphi$ is l.s.c.\
  because $z\in X\setminus F_{n(x)+1}$ implies that $n(z)\leq n(x)$ and, therefore,
  $\varphi(x)=\left[0,2^{-n(x)}\right]\subset
  \left[0,2^{-n(z)}\right]=\varphi(z)$. Finally, observe that $g(x)=1$, $x\in
  X$, is a continuous $1$-selection for $\varphi$. Hence, by (c), $\varphi$ has a
  continuous selection $f:X\to \setR$ such that $|g(x)-f(x)| =1-f(x)<1$ for
  every $x\in X$, or, in other words, $f(x)>0$ for all $x\in X$. Finally,
  define $W_n=f^{-1}\left(\left(-\infty,2^{-n+1}\right)\right)$,
  $n<\omega$. Thus, we get a sequence $\{W_n:n<\omega\}$ of open open subsets of
  $X$ such that $F_n\subset W_n$ for every $n<\omega$. Indeed, $x\in F_n$ implies
  $n\leq n(x)$, so $f(x)\in \varphi(x)=\left[0,2^{-n(x)}\right]\subset
  \left[0,2^{-n}\right]\subset \left[0,2^{-n+1}\right]$. Since
  $\overline{W_n}\subset f^{-1}\left(\left(-\infty,2^{-n+1}\right]\right)$, $n<\omega$,
  and $f(x)>0$ for every $x\in X$, we have that
  $\bigcap\left\{\overline{W_n}:n<\omega\right\}=\emptyset$. That is, $X$ is countably
  paracompact, see \cite[Theorem 5.2.1]{engelking:89}.
\end{proof}

For collectionwise normal spaces we have a very similar result which,
in particular, illustrates the difference with countably paracompact
ones (see, (c) of Theorem
\ref{theorem-countably-paracompact-controlled}).

\begin{proposition}
  \label{proposition-convex-approximation}
  Let $E$ be a Banach space, $X$ be a collectionwise normal space,
  $\varphi:X\to \mathscr{C}'(E)$ be an l.s.c.\ convex-valued mapping, and let
  $g:X\to E$ be a continuous $\varepsilon$-selection for $\varphi$ for some $\varepsilon>0$. Then,
  $\varphi$ has a continuous selection $f:X\to E$ such that $d(f(x),g(x))\leq \varepsilon$
  for every $x\in X$.
\end{proposition}

\begin{proof}
  Define a mapping $\psi:X\to \mathscr{F}(E)$ by
  $\psi(x)=\overline{B_\varepsilon(g(x))}$, $x\in X$. Then, $\psi$ is convex-valued and
  $d$-proximal continuous in the sense of \cite{gutev:95c}. Define
  another mapping $\theta:X\to \mathscr{F}(E)$ by $\theta(x)=\overline{\varphi(x)\cap
    B_\varepsilon(g(x))}$, $x\in X$. According to \cite[Propositions 2.3 and
  2.5]{michael:56a}, $\theta$ is l.s.c.\ and clearly it is also
  convex-valued. Finally, observe that $\theta(x)\subset \psi(x)$ for every $x\in X$,
  while $\theta(x)\neq \psi(x)$ implies that $\theta(x)$ is compact. Then, by
  \cite[Lemma 4.2]{gutev-ohta-yamazaki:00}, $\theta$ has a continuous
  selection $f:X\to E$. This $f$ is as required.
\end{proof}

  Following the idea Proposition
  \ref{proposition-convex-approximation}, one can extend Theorem
  \ref{theorem-countably-paracompact-controlled} to the case of
  countably paracompact and $\tau$-collectionwise normal spaces.

  \begin{theorem}
    \label{theorem-generalized-countably-paracpmpact}
    For a $T_1$-space $X$ and an infinite cardinal number $\tau$, the
    following are equivalent\textup{:}
    \begin{enumerate}
    \item[(a)] $X$ is  countably
    paracompact and $\tau$-collectionwise normal.
  \item[(b)] If $E$ is a Banach space with $w(E)\leq \tau$, $\varphi:X\to
    \mathscr{C}'(E)$ is l.s.c.\ and convex-valued, $\eta:X\to (0,+\infty)$ is
    lower semi-continuous, and $g:X\to Y$ is a continuous $\eta$-selection
    for $\varphi$, then $\varphi$ has a continuous selection $f:X\to E$ such that
    $d(f(x),g(x))<\eta(x)$ for all $x\in X$.
    \end{enumerate} 
  \end{theorem}

  \begin{proof}
    (a)$\ \Rightarrow\ $(b). As in (a) of the proof of Theorem
    \ref{theorem-countably-paracompact-controlled}, there exists a
    continuous function $\alpha:X\to (0,+\infty)$ such that
    $d(g(x),\varphi(x))<\alpha(x)<\eta(x)$ for every $x\in X$. Next, as in the proof
    of Proposition \ref{proposition-convex-approximation}, define a
    ($d$-proximal) continuous $\psi(x)=\overline{B_{\alpha(x)}(g(x))}$, $x\in
    X$, and another l.s.c.\ $\theta:X\to \mathscr{F}(E)$ by
    $\theta(x)=\overline{\varphi(x)\cap B_{\alpha(x)}(g(x))}$, $x\in X$. As in the proof of
    \cite[Lemma 4.2]{gutev-ohta-yamazaki:00}, this $\theta$ has the
    Selection Factorization Property in the sense of
    \cite{nedev:80}. Hence, by \cite[Proposition 4.3]{nedev:80}, $\theta$
    has a continuous selection $f:X\to E$. According to the definition
    of $\theta$, we get that $d(f(x),g(x))<\eta(x)$ for all $x\in X$.\smallskip

    (b)$\ \Rightarrow\ $(a). This implication is based on standard arguments. In
    fact, $X$ will be countably paracompact and normal by Theorem
    \ref{theorem-countably-paracompact-controlled}. To show that $X$
    is also $\tau$-collectionwise normal, let $\mathscr{D}$ be a discrete
    family of closed subsets of $X$, with $|\mathscr{D}|\leq \tau$, and let
    $\ell_1(\mathscr{D})$ be the Banach space of all functions
    $y:\mathscr{D}\to \mathbb{R}$, with $\sum\big\{|y(D)| : D\in
    \mathscr{D}\big\}<\infty $, equipped with the norm $\|y\|=\sum\big\{|y(D)| :
    D\in \mathscr{D}\big\}$. Also, let $\vartheta(D)=0$, $D\in \mathscr{D}$, be the
    origin of $\ell_1(\mathscr{D})$. For every $D\in \mathscr{D}$, consider
    the function $e_D:\mathscr{D}\to \setR$ defined by $e_D(D)=1$ and
    $e_D(T)=0$ for $T\in \mathscr{D}\setminus\{D\}$.  Then, $e_D\in
    \ell_1(\mathscr{D})$, $D\in \mathscr{D}$, and $\|e_D-\vartheta\| =1$ for every
    $D\in \mathscr{D}$. Finally, define an l.s.c.\ mapping $\varphi:X\to
    \mathscr{C}'(\ell_1(\mathscr{D}))$ by $\varphi(x)=\{e_D\}$ if $x\in D$ for some
    $D\in \mathscr{D}$ and $\varphi(x)=\ell_1(\mathscr{D})$ otherwise. Then, $\varphi$
    is convex-valued and $g(x)=\vartheta$, $x\in X$, is a continuous
    $2$-selection for $\varphi$. Since $w(\ell_1(\mathscr{D}))\leq \tau$, by (b), $\varphi$
    has a continuous selection $f:X\to \ell_1(\mathscr{D})$. Then,
    $U_D=f^{-1}(B_1(e_D))$, $D\in \mathscr{D}$, is a pairwise disjoint
    family of open subsets of $X$ such that $D\subset {U}_D$, $D\in
    \mathscr{D}$. Since $X$ is normal, this implies that it is also
    $\tau$-collectionwise normal.
  \end{proof}

  Exactly the same arguments as those in the proof of Theorem
  \ref{theorem-generalized-countably-paracpmpact} show that if $X$ is
  $\tau$-collectionwise normal, $E$ is a Banach space with $w(E)\leq \tau$,
  $\varphi:X\to \mathscr{C}'(E)$ is an l.s.c.\ convex-valued mapping, $\eta:X\to
  (0,+\infty)$ is continuous, and $g:X\to E$ is a continuous $\eta$-selection,
  then $\varphi$ has a continuous selection $f:X\to E$ such that
  $d(f(x),g(x))\leq \eta(x)$ for all $x\in X$. Motivated by this and
  Proposition \ref{proposition-convex-approximation}, we have the
  following natural question.

  \begin{problem}
    \label{question-controlled-collectionwise-normal}
    Let $X$ be a collectionwise normal space, $E$ be a Banach space,
    $\varphi:X\to \mathscr{C}'(E)$ be an l.s.c.\ convex-valued mapping, and
    let $g:X\to E$ be a continuous $\eta$-selection for $\varphi$ for some lower
    semi-continuous function ${\eta:X\to (0,+\infty)}$. Then, does $\varphi$ have a
    continuous selection $f:X\to E$ with $d(f(x),g(x))\leq \eta(x)$ for all
    $x\in X$?
  \end{problem}

  Let us point out that the answer to Question
  \ref{question-controlled-collectionwise-normal} is ``yes'' if so is
  the answer to Question \ref{question-selection-deficiency}. Indeed,
  if ${\eta:X\to (0,+\infty)}$ is lower semi-continuous and $g:X\to E$ is
  continuous, then the mapping $\psi(x)=B_{\eta(x)}(g(x))$, $x\in X$, will have an
  open graph. If $g$ is also an $\eta$-selection for $\varphi:X\to
  \mathscr{C}'(E)$, then $\theta(x)=\overline{\varphi(x)\cap \psi(x)}$, $x\in X$, will
  have a selection $\mathscr{C}(E)$-deficiency. Finally, if $f:X\to E$
  is a continuous selection for $\theta$, then $d(f(x),g(x))\leq \eta(x)$ for all
  $x\in X$.

\newcommand{\noopsort}[1]{} \newcommand{\singleletter}[1]{#1}
\providecommand{\bysame}{\leavevmode\hbox to3em{\hrulefill}\thinspace}
\providecommand{\MR}{\relax\ifhmode\unskip\space\fi MR }
\providecommand{\MRhref}[2]{%
  \href{http://www.ams.org/mathscinet-getitem?mr=#1}{#2}
}
\providecommand{\href}[2]{#2}


\end{document}